\documentclass[12pt]{amsart}
\usepackage{amsmath,amssymb,latexsym}
\def\dcl{\mathrm{dcl}}
\def\acl{\mathrm{acl}}

\def\bdd{\mathrm{bdd}}

\def\St{\mathrm{stab}}
\def\p{\varphi}
\def\proves{\vdash}
\def\M{\mathfrak M}
\def\N{\mathfrak N}

\def\H{\mathfrak H}
\def\tp{\mathrm{tp}}

\def\Ind#1#2{#1\setbox0=\hbox{$#1x$}\kern\wd0\hbox to 0pt{\hss$#1\mid$\hss}
\lower.9\ht0\hbox to 0pt{\hss$#1\smile$\hss}\kern\wd0}
\def\ind{\mathop{\mathpalette\Ind{}}}
\def\Notind#1#2{#1\setbox0=\hbox{$#1x$}\kern\wd0\hbox to 0pt{\mathchardef
\nn="3236\hss$#1\nn$\kern1.4\wd0\hss}\hbox to 0pt{\hss$#1\mid$\hss}\lower.9\ht0
\hbox to 0pt{\hss$#1\smile$\hss}\kern\wd0}

\theoremstyle{plain}
\newtheorem{theorem}{Theorem}[section]
\newtheorem{proposition}[theorem]{Proposition}
\newtheorem{fact}[theorem]{Fact}
\newtheorem{lemma}[theorem]{Lemma}
\newtheorem{corollary}[theorem]{Corollary}
\newtheorem*{claim}{Claim}

\theoremstyle{definition}
\newtheorem{definition}[theorem]{Definition}
\newtheorem{remark}[theorem]{Remark}
\newtheorem{expl}[theorem]{Example}
\newtheorem{frage}[theorem]{Question}

\def\bsp{\begin{expl}}
\def\ebsp{\end{expl}}
\def\beh{\begin{claim}}
\def\ebeh{\end{claim}}
\def\defn{\begin{definition}}
\def\edefn{\end{definition}}
\def\satz{\begin{theorem}}
\def\esatz{\end{theorem}}
\def\tats{\begin{fact}}
\def\etats{\end{fact}}
\def\kor{\begin{corollary}}
\def\ekor{\end{corollary}}
\def\lmm{\begin{lemma}}
\def\elmm{\end{lemma}}
\def\bem{\begin{remark}}
\def\ebem{\end{remark}}
\def\bew{\par\noindent{\em Proof: }}

\def\satzli{\begin{proposition}}
\def\esatzli{\end{proposition}}

\begin{document}
\title{The right angle to look at \hbox{orthogonal sets}}
\author{Frank O. Wagner}
\address{Universit\'e de Lyon; CNRS; Universit\'e Claude Bernard Lyon 1; Institut Camille Jordan UMR5208, 43 bd du 11 novembre 1918, 69622 Villeurbanne Cedex, France}
\email{wagner@math.univ-lyon1.fr}
\keywords{orthogonality; simplicity; internality}
\date{3 December 2014}
\subjclass[2000]{03C45}
\date{}
\thanks{Partially supported by ANR-09-BLAN-0047 Modig and ANR-13-BS01-0006
ValCoMo}
\begin{abstract}If $X$ and $Y$ are orthogonal hyperdefinable sets such that $X$
is simple, then any group $G$ interpretable in $X\cup Y$ has a
normal hyperdefinable $X$-internal subgroup $N$ such that $G/N$ is
$Y$-internal; $N$ is unique up to commensurability. In order to
make sense of this statement, local simplicity theory for
hyperdefinable sets is developed. Moreover, a version of Schlichting's Theorem for hyperdefinable families of commensurable subgroups is shown.\end{abstract}
\maketitle

\section*{Introduction}
Two definable sets $X$ and $Y$ in some structure are said to be {\em orthogonal}
if every definable subset of $X\times Y$ is a finite union of {\em rectangles},
i.e.\ of subsets of the form $U\times V$ with $U\subseteq X$ and $V\subseteq Y$
definable. It follows that if $X$ and $Y$ are orthogonal groups, every
definable subgroup $H$ of $X\times Y$ has a subgroup of finite index of the form 
$U\times V$ with $U\le X$ and $V\le Y$ subgroups: As $H$ is a finite union of rectangles, one can find a maximal definable rectangle $U\times V\subseteq H$ containing the identity $1=(1_X,1_Y)$. As $H$ also contains
$$(U\times V)^{-1}(U\times V)=(U^{-1}\times V^{-1})(U\times V)=U^{-1}U\times V^{-1}V\supseteq U\times V,$$
we obtain $U^{-1}U=U$ and $V^{-1}V=V$ by maximality, so $U\le X$ and $V\le Y$ are subgroups; moreover $U\times V$ is unique. Any other maximal rectangle contained in $H$ can be translated to contain $1$, and must thus be a coset of $U\times V$. So $U\times V$ has finite index.

However, the situation is considerably more complicated for
a group $G$ definable, or more generally interpretable, in $X\cup Y$, as it need
not be a direct product of a group interpretable in $X$ and a group interpretable in $Y$.
In fact, an example by Berarducci and Mamino \cite[Example 1.2]{BM} shows that
$G$ need not have {\em any} subgroup interpretable in either $X$ or $Y$.
However, they prove \cite[Theorem 7.1]{BM} that if $X$ is superstable of finite
and definable Lascar rank, then any group $G$ interpretable in $X\cup Y$ has a
normal subgroup $N$ interpretable in $X$, such that $G/N$ is interpretable in
$Y$.

In this paper we shall generalize their result to the case where $X$ is merely
simple. In this context, definability has to be replaced by type-definability,
as even for a definable group the tools of simplicity theory in general only
yield type-definable subgroups. In fact, we even have to study {\em
hyperdefinable} groups, since the quotient $G/N$, for $N$ type-definable, will
be of that form. We therefore put ourselves in the hyperdefinable context and
assume right from the start that our orthogonal sets $X$ and $Y$ are merely
hyperdefinable. To this end, we shall include a quick development of 
hyperdefinability in section \ref{sec1}, and of local simplicity theory for 
hyperdefinable sets in section \ref{sec4}. Moreover, the general theory only yields $N$ unique up to commensurability. Since there is {\em a priori} no simple hyperdefinable set containing all conjugates of $N$, we cannot use the usual locally connected
component from simplicity theory \cite[Definition 4.5.15]{wa00}. We therefore show a completely general version of Schlichting's Theorem for a hyperdefinable family of commensurable subgroups in section \ref{secfive}. Note that we do recover the theorem
by Berarducci and Mamino even for general supersimple (and definable) $X$
(Corollary \ref{BM}).

Another problem is that of parameters. The usual hypothesis would be that of
{\em stable embedding}, i.e.\ that every hyperdefinable subset of $X$ is
hyperdefinable with parameters in $X$. We shall circumvent this issue by only
ever considering parameters from $X\cup Y$, as orthogonality automatically
yields stable embeddedness of $X$ and of $Y$ in $X\cup Y$. 

We shall work in a big $\kappa$-saturated and strongly $\kappa$-homogeneous
monster model $\M$, where $\kappa$ is bigger than any cardinality we wish to
consider. We shall not usually distinguish between elements and tuples.

{\em Acknowledgements.} I should like to thank the anonymous referee whose comments have greatly contributed to improve the presentation of the paper, and induced me to include Section \ref{secfive} on Schlichting's Theorem (and to prove the results therein).

\section{Hyperimaginaries}\label{sec1}
\defn\label{defnHI}\begin{itemize}\item A {\em countable equivalence relation} is an 
equivalence relation given by the conjunction of countably many formulas (and 
hence only using countably many parameters and variables); it is {\em over} a 
set $A$ of parameters if the formulas only use parameters from $A$. 
\item A {\em hyperimaginary} (element) of {\em type $E$} is the class $a_E$ of 
some 
tuple $a$ (of the right length) modulo a countable equivalence relation $E$ over 
$\emptyset$. 
\item A hyperimaginary $e$ is {\em definable} over some set $B$ of 
hyperimaginaries if every automorphisms of the monster model which fixes $B$ 
pointwise fixes also $e$; it is {\em bounded} over $B$ if its orbit under the 
group of automorphisms fixing $B$ pointwise has bounded size (smaller than the 
saturation degree of the monster model). The {\em hyperimaginary definable 
closure} $\dcl^{heq}(B)$ of $B$ is the set of hyperimaginaries definable over 
$B$; the {\em hyperimaginary bounded closure} $\bdd(B)$ of $B$ is the set of 
all hyperimaginaries bounded over $B$. Clearly, both $\dcl^{heq}$ and 
$\bdd$ are idempotent operators.\footnote{By compactness, an imaginary element 
in $\bdd(B)$ is already in the algebraic closure $\acl(B)$. So there is no need 
for a superscript $\bdd^{heq}$.}
\item Two hyperimaginaries are {\em equivalent} if they are interdefinable.
\item If $e$ is hyperimaginary, a {\em representative} for $e$ is any real (or
imaginary) tuple $a$ with $e\in\dcl^{heq}(a)$.\end{itemize}
\edefn

\bem\label{eqnrel}\begin{enumerate}\item If $E$ is an arbitrary type-definable 
equivalence relation over $\emptyset$ (given by an intersection of arbitrary 
size, on tuples of arbitrary 
length), it is easy to see that $E$ is equivalent to a conjunction of
subintersections $E_i$, each one defining a countable equivalence 
relation on a countable subtuple $x_i$. So 
$$\qquad\quad aEb\qquad\Leftrightarrow\qquad\forall i\ a_iE_ib_i.$$
This means that an automorphism fixes $a_E$ if and only iff it fixes 
$(a_i)_{E_i}$ for all $i$, and we can replace $a_E$ by the sequence 
$\big((a_i)_{E_i}\big)_i$ of hyperimaginaries.
\item If $E_A$ is a countable equivalence relation over $A$, we 
consider the countable equivalence relation
$$\qquad\quad\,xx'Fyy'\quad\Leftrightarrow\quad 
\big(x'=y'\land x'\models\tp(A)\land xE_{x'}y\big)\lor xx'=yy'$$
where we only consider the countable subset of $A$ actually occurring in the 
definition of $E_A$, and $E_{x'}$ is the result of substituting $x'$ for $A$ in the definition of $E_A$. Then for any $a$ the class $a_{E_A}$ is fixed by an 
automorphism fixing $A$ iff and only iff $(aA)_F$ is fixed, and we can use 
the hyperimaginary $(aA)_F$ instead of $a_{E_A}$.
\item If $E$ is a countable partial type over $A$ which defines an 
equivalence relation on some partial type $\pi$ over $A$, then by compactness 
there is a countable subtype $\pi_0$ such that $E$ defines an equivalence 
relation on $\pi_0$. Then
$$\qquad\quad\ xFy\quad\Leftrightarrow\quad\big(\pi_0(x)\land\pi_0(y)\land 
xEy\big)\lor x=y$$
is a countable equivalence relation over $A$ extending $E$.
\end{enumerate}\ebem

If $E=(E_i:i\in I)$ is a sequence of countable equivalence 
relations over $\emptyset$ and $a=(a_i:i\in I)$ is a
sequence of tuples of the right length, we put $a_E=((a_i)_{E_i}:i\in I)$, and we say that $a_E$ is a tuple of hyperimaginaries. 
Similarly, we write $aEb$ if $a_iE_ib_i$ for all $i\in I$. Note that Remark \ref{eqnrel} justifies that we restrict to countable equivalence relations over $\emptyset$ in Definition \ref{defnHI}: Indeed, any other equivalence class one might wish to consider is just a tuple of hypermaginaries.

\defn Let $a_E$ and $b_F$ be tuples of hyperimaginaries. The {\em type}
$\tp(a_E/b_F)$ is given by all partial types over $b$ of the form 
$$\exists yz\ [xEy\land zFb\land\p(y,z)]$$
true of $a$, where $\p$ is a parameter-free formula. It is easy to see that (in the monster model) two 
tuples of hyperimaginaries of type $E$ have the same type over $b_F$ if and 
only if they are conjugate by an automorphism fixing $b_F$.\edefn

Note that the type of a hyperimaginary over $b_F$ is just a maximal
$E$-invariant partial real type over $b$ invariant under automorphisms 
fixing $b_F$. For any two representatives of $b_F$, any two 
such types are equivalent. We shall say that a partial type $\pi(y)$
is a partial $E$-type if $\pi(y)$ is $E$-invariant.

\defn A set $X$ is {\em hyperdefinable} over some parameters $A$ if it is of the
form $Y/E$, where $Y$ is a type-definable set in countably many
variables and $E$ a countable equivalence relation on $Y$, both over $A$.
We denote by $X_A^{heq}$ the collection of all hyperimaginaries in the 
definable closure of $A$ and some tuple from $X$. If $A=\emptyset$ it is omitted.\edefn

For the rest of the paper, all tuples and parameter sets are 
hyperimaginary, unless stated otherwise. We shall not distinguish between
elements and tuples of elements from a set.

\section{Orthogonality}
\defn Let $X$, $Y$ be $A$-hyperdefinable sets in some structure $\M$. We
say that $X$ and $Y$ are {\em orthogonal} over $A$, denoted $X\perp_A Y$, if for
any tuples $a$ from $X$ and $b$ from $Y$, the partial type
$\tp(a/A)\cup\tp(b/A)$ determines $\tp(ab/A)$. If $A=\emptyset$ it will be omitted.\edefn
\bem Note that we do not require $X$ (or $Y$) to be {\em stably embedded}, 
i.e.\ that every hyperdefinable subset of $X$ be hyperdefinable with 
parameters in $A\cup X$. In a stable theory, every hyperdefinable subset is 
stably embedded, but this need not hold in general. We shall compensate for 
the lack of stable embeddedness by restricting our additional parameters to 
$X_A^{heq}\cup Y_A^{heq}$.\footnote{See Proposition \ref{perp}. In a stable 
theory, $X\perp_A Y$ implies $X\perp_B Y$ for any $B\supseteq A$ (full 
orthogonality).}\ebem
\bsp If $\M_1$ and $\M_2$ are two structures and $\N=\M_1\times\M_2$ with a
predicate $X$ for $\M$ and a predicate $Y$ for $\N$, then $X$ and $Y$ are
orthogonal in $\N$ over $\emptyset$ (and in fact over any set of parameters).\ebsp
\bem If $X$ and $Y$ are orthogonal type-definable sets over $A$ and $Z\subseteq 
X^k\times Y^\ell$ is relatively $A$-definable, then $Z$ is a finite union of 
rectangles $A_i\times B_i$, where $A_i\subseteq X^k$ and $B_i\subseteq Y^\ell$
are relatively $A$-definable.\ebem
\bew For any $z=(x,y)\in Z$ we have that 
$$\tp(x/A)\cup\tp(y/A)\proves (x,y)\in Z.$$
By compactness there are relatively $A$-definable subsets $A_z\subseteq X^k$ in
$\tp(x/A)$ and $B_z\subseteq Y^\ell$ in $\tp(y/A)$ with $A_z\times B_z\subseteq
Z$. 
Again by compactness, finitely many of these rectangles suffice to cover 
$Z$.\qed
\bem\label{bdd} $X\perp_A X$ if and only if 
$X\subseteq\dcl^{heq}(A)$.\ebem
\bew For any $x,x'\in X$ we have
$$\tp(x/A)\cup\tp(x'/A)\vdash\tp(x,x'/A).$$
If $x\notin\dcl^{heq}(A)$ choose $x'\equiv_A x$ with $x'\not=x$. Then
$xx'\equiv_A xx$, a contradiction.\qed

For the rest of this section, $X$ and $Y$ will be orthogonal
$\emptyset$-hyper\-definable sets. We note first that orthogonality is 
preserved under adding parameters from $X^{heq}\cup Y^{heq}$, and 
interpretation:
\satzli\label{perp} If $X'\subseteq X^{heq}$ and $Y'\subseteq Y^{heq}$ are
hyperdefinable
over some parameters $A\subseteq X^{heq}\cup Y^{heq}$, then $X'\perp_A
Y'$.\esatzli 
\bew Suppose $A=(a,b)$ with $a\in X^{heq}$ and $b\in Y^{heq}$, and consider
tuples $a'\in X'$ and $b'\in Y'$. Choose representatives $\bar a,\bar a'\in X$
of $a,a'$ and $\bar b,\bar b'\in Y$ of $b,b'$. Then $\tp(\bar a\bar
a')\cup\tp(\bar b\bar
b')\proves\tp(\bar a\bar a'\bar b\bar b')$.

Now if $a''\equiv_Aa'$ and $b''\equiv_Ab'$, we can find $A$-conjugates $\tilde
a\bar a''$ of $\bar a\bar a'$ and $\tilde b\bar b''$ of $\bar b\bar b'$ such
that $a''\tilde a\bar a''\equiv_A a'\bar a\bar a'$ and $b''\tilde b\bar
b''\equiv_A b'\bar b\bar b'$. By orthogonality of $X$ and $Y$, we obtain $\bar
a''\tilde a\bar b''\tilde b\equiv \bar a'\bar a\bar b'\bar b$, whence
$a''ab''b\equiv a'ab'b$, and thus $a''b''\equiv_A a'b'$.\qed

\satzli\label{stabemb} $X$ is stably embedded in $X\cup Y$: For tuples $a\in
X^{heq}$ and $b\in Y^{heq}$, every $ab$-hyperdefinable
subset $X'$ of $X^{heq}$ is hyperdefinable over $a$.\esatzli 
\bew If $\Phi(x,a,b)$ hyperdefines $X'$ and $\Psi(y)=\tp(b)$, put
$$\Phi'(x,a)=\exists y\,[\Psi(y)\land\Phi(x,a,y)].$$
Clearly $\Phi(x,a,b)\proves\Phi'(x,a)$. Conversely, suppose
$a'\models\Phi'(x,a)$, and choose $b'\models\Psi$ with
$a'\models\Phi(x,a,b')$. By orthogonality $a'ab\equiv a'ab'$, whence
$a'\models\Phi(x,a,b)$, and
$\Phi'(x,a)$ hyperdefines $X'$.\qed

We put $\dcl^{heq}_X(A)=\dcl^{heq}(A)\cap X^{heq}$ and $\bdd_X(A)=\bdd(A)\cap
X^{heq}$.
\kor\label{dclbdd} Suppose $a\in X^{heq}$ and $b\in Y^{heq}$. Then
$$\dcl^{heq}_X(a,b)=\dcl^{heq}_X(a)\quad\text{and}\quad
\bdd_X(a,b)=\bdd_X(a).$$\ekor 
\bew Immediate from Proposition \ref{stabemb}: If $X'$ is a singleton 
(resp.\ bounded) subset of $X^{heq}$ hyperdefinable over $ab$ containing some 
element $e\in\dcl^{heq}_X(ab)$ (resp.\ $e\in\bdd_X(ab)$), then $X'$ is 
hyperdefinable already over~$a$.\qed

\section{Weak elimination of hyperimaginaries}
In this section, $X$ and $Y$ will be $\emptyset$-hyperdefinable sets.
\defn Let $Z$ be $\emptyset$-hyperdefinable. We say that $Z$ has
{\em weak elimination of hyperimaginaries with respect to $X^{heq}$ and
$Y^{heq}$} if for every $z\in Z^{heq}$ there is some $x\in\bdd_X(z)$ and
$y\in\bdd_Y(z)$ with $z\in\dcl^{heq}(xy)$.\edefn
For the rest of this section, $X$ and $Y$ will be orthogonal over $\emptyset$.
\satz\label{wehi} The set $X\cup Y$ has weak
elimination of hyperimaginaries with respect to $X^{heq}$ and $Y^{heq}$.\esatz
\bew Consider $z\in(X\cup Y)^{heq}$, say $z=(x,y)_E$ for some tuples 
$x\in X$, $y\in Y$ and countable equivalence relation $E$ over $\emptyset$. 
For $x'\equiv x$ consider the hyperdefinable equivalence relation $E_{x'}$ on 
$\tp(y)$ given by 
$$yE_{x'}y'\quad\Leftrightarrow\quad (x',y)E(x',y').$$
Then $E_{x'}$ is $\emptyset$-hyperdefinable by Proposition \ref{stabemb}, and 
does not depend on the choice of $x'\equiv x$. Similarly, for $y'\equiv y$ the 
equivalence relation 
$$xE_{y'}x'\quad\Leftrightarrow\quad (x,y')E(x',y')$$
on $\tp(x)$ does not depend on $y'\equiv y$. Clearly 
$z\in\dcl^{heq}(x_{E_{y'}},y_{E_{x'}})$.

We claim that $x_{E_{y'}}$ is bounded over $z$.
If not, there is an indiscernible sequence
$(x_i,y_i:i<\omega)$ in $\tp(x,y/z)$ with $\neg x_i E_{y'}x_j$ for $i\not=j$. By 
orthogonality, for $i<j$,
$$\tp(x_i,x_j)\cup\tp(y_i,y_j)\proves\tp((x_i,y_i),(x_j,
y_j)).$$
But $\tp(x_i,x_j)=\tp(x_i,x_k)$ for $i<k<j$, whence
$$\tp((x_i,y_i),(x_j,y_j))=\tp((x_i,y_i),(x_k,y_j)).$$
Now $(x_i,y_i)E(x,y)E(x_j,y_j)$ holds since $z=(x,y)_E$. Hence
$$(x_k,y_j)E(x_i,y_i)E(x_j,y_j).$$
Thus $x_kE_{y'}x_j$, a contradiction.

Hence $x_{E_{y'}}\in\bdd_X(z)$; similarly 
$y_{E_{x'}}\in\bdd_Y(z)$.\qed

\kor\label{XY=X+Y} For any set $A$ of parameters, $\bdd_{XY}(A)$ and 
$\bdd_X(A)\cup\bdd_Y(A)$ are interdefinable. Moreover, for $aA\subset (X\cup 
Y)^{heq}$ we have $\tp(a/\bdd_{XY}(A))\proves\tp(a/\bdd(A))$.\ekor
\bew Clearly $\bdd_X(A)\cup\bdd_Y(A)\subseteq\bdd_{XY}(A)$.

For the converse inclusion, let $z\in\bdd_{XY}(A)$. By 
Theorem \ref{wehi} there is $x\in\bdd_X(z)$ and $y\in\bdd_Y(z)$ with 
$z\in\dcl^{eq}(xy)$. So
$$z\in\dcl^{heq}(\bdd_X(\bdd(A)),\bdd_Y(\bdd(A)))=\dcl^{heq}(\bdd_X(A),
\bdd_Y(A)).$$

For the second assertion, let $B$ be a set of representatives for $\bdd(A)$ and $F$ a type-definable equivalence relation such that $B_F$ is equivalent to $\bdd(A)$. Then equality of $E$-type over $\bdd(A)$ is a bounded equivalence relation $E_B$  type-definable over $B$, given by
$$\begin{aligned}xE_By\quad\Leftrightarrow\bigwedge_{\p\text{ a 
$B$-formula}}&\big[\big(\p(x,B)\to\exists y'z\ [yEy'\land 
BFz\land\p(y',z)]\big)\\
&\land\big(\p(y,B)\to\exists x'z\ [xEx'\land 
BFz\land\p(x',z)]\big)\big].\end{aligned}$$
As $E_B$ is invariant under any $A$-automorphism, it is in fact type-definable 
over $A$. By Remark \ref{eqnrel} the class $a_{E_B}$ is interdefinable with a 
tuple
$$((aA_i)_{E_i}:i\in I)\in\bdd_{XY}(A),$$
where the $A_i\subseteq A$ are countable. 
Since the partial type $(xA_i)E_i(aA_i)$ is in $\tp(a/\bdd_{XY}(A))$ for all 
$i\in I$, we get the result.\qed

\kor\label{bddperp}If $X'\subset\bdd(X)$ and $Y'\subset\bdd(Y)$
are $\emptyset$-hyperdefinable, then $X'\perp_{\bdd_{XY}(\emptyset)} Y'$.\ekor
\bew It is clearly sufficient to show $X\perp_{\bdd_{XY}(\emptyset)} Y'$. Given
$x\in X$ and $y\in Y'$, consider $y_0\in Y$ with $y\in\bdd(y_0)$, and put $\bar
y=\bdd_Y(y_0)$. Now if 
$$x'\equiv_{\bdd_{XY}(\emptyset)}x\qquad\text{and}\qquad y'\equiv_{\bdd_{XY}
(\emptyset)}y,$$
choose $\bar y'$ with $\bar y'y'\equiv_{\bdd_{XY}(\emptyset)}\bar yy$. As
$X\perp_{\bdd_{XY}(\emptyset)}Y$ by Lemma \ref{perp}, we have $x\bar
y\equiv_{\bdd_{XY}(\emptyset)}x'\bar y'$. Since
$$\bdd_{XY}(\bar y)\in\dcl^{heq}(\bdd_X(\bar y),\bdd_Y(\bar y))
=\dcl^{heq}(\bdd_X(\emptyset),\bar y)$$
by Corollary \ref{dclbdd}, we obtain
$$x\bdd_{XY}(\bar y)\equiv x'\bdd_{XY}(\bar y')\quad\text{and}\quad 
\bdd_{XY}(\bar y)y\equiv \bdd_{XY}(\bar y')y'.$$
(Note that $\bdd_{XY}(\emptyset)$ is part of the tuples on either side, so we do not have to work over it.) Choose $x''$ with $x\bdd_{XY}(\bar y)y\equiv x''\bdd_{XY}(\bar y')y'$. Then
$$x''\bdd_{XY}(\bar y')\equiv x\bdd_{XY}(\bar y)\equiv x'\bdd_{XY}(\bar
y'),$$
so $\tp(x''/\bdd_{XY}(\bar y'))=\tp(x'/\bdd_{XY}(\bar y'))$. By Corollary
\ref{XY=X+Y} we obtain $\tp(x''/\bdd(\bar y'))=\tp(x'/\bdd(\bar y'))$. As
$y\in\bdd(\bar y)$, we get in particular
$$x'y'\equiv x''y'\equiv xy.$$
The result follows.\qed

\bsp\label{bsp} We do need $\bdd_{XY}(\emptyset)$ in Corollary \ref{bddperp}, 
as we might take $X'=X\times\bdd_Y(\emptyset)$. Then $X'\not\perp Y$ unless 
$\bdd_Y(\emptyset)=\dcl^{heq}_Y(\emptyset)$.\ebsp

\section{Internality and analysability}
\defn Let $X$ and $Y$ be hyperdefinable sets over $A$. We say that $X$ is
{\em (almost) $Y$-internal} if there is some parameter set $B$ such that for
every $a\in X$ there is a tuple $b\in Y$ with $a\in\dcl^{heq}(Bb)$
(or $a\in\bdd(Bb)$, respectively).\footnote{In simplicity theory, this is
called {\em finite generation\/}; for {\em internality} we would require for
every $a\in X$ the existence of some $B\ind_A a$ and tuple $b\in Y$ with
$a\in\dcl^{heq}(Bb)$.}\par
If the parameters $B$ can be chosen in some set $Z$, we say that
$X$ is (almost) $Y$-internal {\em within $Z$}.\par
We say that $X$ is {\em $Y$-analysable (within $Z$)} if for all $a\in X$ there 
is a sequence $(a_i:i<\alpha)$ such that $\tp(a_i/A,a_j:j<i)$ is 
$Y$-internal (within $Z$)
for every $i<\alpha$, and $a\in\bdd(A,a_i:i<\alpha)$.\edefn
For the rest of the section, $X$ and $Y$ will be hyperdefinable orthogonal sets 
over $\emptyset$.
\satzli\label{anaint} If an $\emptyset$-hyperdefinable set $X'$ is 
$X$-analysable within $X\cup Y$, then $X'$ is almost $X$-internal within
$\bdd_Y(\emptyset)$; if $X'$ is $X$-internal within $X\cup Y$, then $X'$ is
$X$-internal within $\bdd_Y(\emptyset)$.\esatzli
\bew We first show that if $X'$ is (almost) $X$-internal within $X\cup Y$, then it is (almost) $X$-internal within $\bdd_Y(\emptyset)$.
So suppose $\bar a\in X$ and $\bar b\in Y$ are such that for every $c_E\in
X'$ there is a tuple $a\in X$ with $c_E\in\bdd(\bar a\bar ba)$. Let 
$\Phi(x,\bar a\bar ba)$ be the $E$-type $\tp(c_E/\bar a\bar 
ba)$. Then for every symmetric formula $\psi(x,y)\in E$ there is 
$n_\psi<\omega$ such that a maximal $\psi$-antichain in $\Phi$ has size 
$n_\psi$, and a formula $\phi_\psi(x,\bar a\bar ba)\in\Phi$ such that every 
$\psi$-antichain in $\phi_\psi$ has size at most $n_\psi$. Consider the 
type-definable relation $F$ on $\tp(\bar a\bar ba)$ given by
$$\begin{aligned}(\bar a'\bar b'a')F(\bar a''\bar 
b''a'')&\Leftrightarrow\!\!\bigwedge_{\psi\in E}\!\!\big[\forall x
\big(\phi_\psi(x,\bar a'\bar b'a')\to\exists x'[\Phi(x',\bar 
a''\bar b''b'')\land x\psi^2x']\big)\\
&\quad\land\forall x'\big(\phi_\psi(x',\bar a''\bar b''a'')\to\exists x
[\Phi(x,\bar a'\bar b'b')\land x\psi^2x']\big)\big],\end{aligned}$$
where $x\psi^2x'$ means $\exists x''\,[\psi(x,x'')\land\psi(x'',x')]$.
Then $(\bar a'\bar b'a')F(\bar a''\bar b''a'')$ holds if and only if 
$\Phi(x,\bar a'\bar b'a')$ and $\Phi(x,\bar a''\bar b''a'')$ contain the same 
points modulo $E$: If they contain the same points modulo $E$, for every 
$\psi$ in $E$ let $(x_i:i<n_\psi)$ be a $\psi$-antichain in $\Phi(x,\bar a'\bar 
b'a')$, and choose $(x'_i:i<n_\psi)$ in $\Phi(x,\bar a''\bar 
b''a'')$ with $x_iEx'_i$ for all $i<n_\psi$. Then whenever $x$ satisfies 
$\phi_\psi(x,\bar a'\bar b'a')$ there is $i<n_\psi$ with $\psi(x,x_i)$, 
whence $x\psi^2x'_i$. By symmetry the converse also holds, so $(\bar a'\bar 
b'a')F(\bar a''\bar b''a'')$. On the other hand, if there is $x$ such that 
$\Phi(x,\bar a'\bar b'a')$ but $\neg xEx'$ for all $x'$ with 
$\Phi(x',\bar a''\bar b''a'')$, by compactness there is $\psi'\in E$ such that
$\neg\psi'(x,x')$ for all $x'$ with 
$\Phi(x',\bar a''\bar b''a'')$. Then any $\psi\in E$ with 
$\psi^2\vdash\psi'$ witnesses $\neg(\bar a'\bar b'a')F(\bar a''\bar 
b''a'')$.

It follows that $F$ is an equivalence relation, and any automorphism fixes 
$(\bar a\bar ba)_F$ if and only if it permutes the set $C\subset X'$ of 
$E$-classes in $\Phi(x,\bar a\bar ba)$. In particular $(\bar 
a\bar ba)_F\in\dcl^{heq}(C)$ and $C\subseteq\bdd((\bar a\bar ba)_F)$, 
whence in particular $c_E\in\bdd((\bar a\bar ba)_F)$. Moreover, if $X'$ 
is 
$X$-internal, then $C=\{c_E\}$ and $c_E\in\dcl^{heq}(\bar a\bar ba)_F)$. 

By weak elimination of hyperimaginaries, there is $\tilde a\in\bdd_X((\bar 
a\bar ba)_F)$ and $\tilde b\in\bdd_Y((\bar 
a\bar ba)_F)$ with $(\bar a\bar ba)_F\in\dcl^{heq}(\tilde a\tilde b)$. Thus we 
are done if we can show $\tilde b\in\bdd_Y(\emptyset)$.

Suppose $\tilde b\notin\bdd_Y(\emptyset)$. Then there is an
$\emptyset$-conjugate $\tilde b'$ of $\tilde b$ outside $\bdd(\bar a\bar b)$; if
$\sigma$ is an automorphism mapping $\tilde b'$ to $\tilde b$, put $\bar
a'\bar b'=\sigma(\bar a\bar b)$. Then $\tilde b\notin\bdd(\bar a'\bar b')$. On
the other hand, since $\bar a'\bar b'\equiv\bar a\bar b$, for every $e\in 
C$ there is $a_e\in X$ with $e\in\bdd(\bar a'\bar b'a_e)$. Therefore
$$\tilde b\in\bdd_Y((\bar a\bar ba)_F)\subseteq\bdd_Y(C)\subseteq\bdd_Y(\bar 
a',\bar b',a_e:e\in C),$$
whence $\tilde b\in\bdd_Y(\bar b')$ by Corollary \ref{dclbdd}, a
contradiction.

Now assume that $x\in X'$ and $(x_i:i<\alpha)$ is an $X$-analysis of $x$ within
$X\cup Y$. We show inductively on $i\le\alpha$ that $\tp(x_j:j<i)$ is
$X$-internal within
$\bdd_Y(\emptyset)$. So suppose $\tp(x_j:j<k)$ is $X$-internal within
$\bdd_Y(\emptyset)$ for all $k<i$. If $i$ is limit, then clearly $\tp(x_j:j<i)$
is $X$-internal within $\bdd_Y(\emptyset)$. If $i=k+1$, then by the result for 
internality $\tp(x_k/x_j:j<k)$ is $X$-internal within $\bdd_Y(x_j:j<k)$ and 
there is $a\in X$ with $$x_k\in\dcl^{heq}(a,\bdd_Y(x_j:j<k),x_j:j<k).$$
Or, by $X$-internality of
$\tp(x_j:j<k)$ within $\bdd_Y(\emptyset)$ there is $a'\in X$ with
$(x_j:j<k)\in\dcl^{heq}(a',\bdd_Y(\emptyset))$. Then by Corollary \ref{dclbdd}
$$\bdd_Y(x_j:j<k)\subseteq\bdd_Y(a',\bdd_Y(\emptyset))=\bdd_Y(\emptyset),$$
and $x_k\in\bdd(a,a',\bdd_Y(\emptyset))$. So
$\tp(x_j:j<i)$ is $X$-internal within $\bdd_Y(\emptyset)$, and $\tp(x)$ is
almost $X$-internal within $\bdd_Y(\emptyset)$.\qed

\kor\label{perpint} Let $X'$ and $Y'$ be $\emptyset$-hyperdefinable. If $X'$ is 
almost $X$-internal within $X\cup Y$ and $Y'$ is almost $Y$-internal within
$X\cup Y$, then $X'\perp_{\bdd_{XY}(\emptyset)} Y'$.\ekor
\bew Proposition \ref{perp} and Corollary \ref{XY=X+Y} yield 
$X\perp_{\bdd_{XY}(\emptyset)}Y$. By Pro\-po\-sition \ref{anaint} we have 
$$X'\subset\bdd(X,\bdd_{XY}(\emptyset))\qquad\text{and}\qquad
Y'\subset\bdd(Y,\bdd_{XY}(\emptyset)).$$
Hence 
$X'\perp_{\bdd_{XY}(\emptyset)}Y'$ by Corollary \ref{bddperp}.\qed

\kor\label{bddint} If an $\emptyset$-hyperdefinable set $Z$ is almost
$X$- and almost $Y$-internal within $X\cup Y$, then it is bounded.\ekor
\bew We have $Z\perp_{\bdd_{XY}(\emptyset)}Z$ by
Corollary \ref{perpint}, so $Z$ is bounded by Remark \ref{bdd}.\qed

\kor\label{aint=>int} If
$Z\subseteq(X\cup Y)^{heq}$ is $\emptyset$-hyperdefinable and almost
$X$-internal within $X\cup Y$, then it is $X$-internal within
$\bdd_Y(\emptyset)$.\ekor
\bew Let $z\in Z$. By weak elimination of hyperimaginaries there is 
$x\in\bdd_X(z)$ and $y\in\bdd_Y(z)$ with $z\in\dcl^{heq}(xy)$. Then
$\tp(y)$ is $Y$-internal since $y\in Y^{heq}$, but also almost $X$-internal, 
as $y\in\bdd(z)$ and $\tp(z)$ is almost $X$-internal. So 
$y\in\bdd_Y(\emptyset)$ by Corollary \ref{bddint}.\qed

Again Example \ref{bsp} shows that we need $\bdd_Y(\emptyset)$ in Corollaries 
\ref{perpint} and \ref{aint=>int}.

\section{Local simplicity}\label{sec4}
\defn Let $A\subseteq B$, and $\pi(x,B)$ be a partial type over $B$. We say that
$\pi(x,B)$ does not divide over $A$ if for any indiscernible sequence
$(B_i:i<\omega)$ in $\tp(B/A)$ the partial type 
$$\bigcup_{i<\omega}\pi(x,B_i)$$
is consistent. Clearly, $\tp(a/B)$ divides over $A$ if and only if
$\tp(a_0/B)$ does so for some finite subtuple $a_0\subseteq a$.\edefn
\bsp\label{ind} If $\tp(a)\perp\tp(b)$, then $\tp(a/b)$ does not divide over
$\emptyset$.\ebsp
We now define the appropriate version of local rank. We follow
Ben Yaacov's
terminology \cite[Definition 1.4]{BY}, more general than \cite[Definition
4.3.5]{wa00}. 
\defn\label{Ddef} Let $\pi(x)$, $\Phi(x,y)$ and $\Psi(y_1,\ldots,y_k)$ be
partial types in (at most) countably many variables.\begin{enumerate} 
\item $\Psi$ is a {\em $k$-inconsistency witness} for $\Phi$ if 
$$\models\forall y_1\ldots y_k\,\neg\exists
x\,[\Psi(y_1,\ldots,y_k)\land\bigwedge_{i=1}^k\Phi(x,y_i)].$$
\item\label{Ddefdef} Let $\Psi$ be a $k$-inconsistency witness for $\Phi$. The 
{\em local $(\Phi,\Psi)$-rank} $D(.,\Phi,\Psi)$ is
defined on partial types in $x$ as follows:\begin{itemize} 
\item $D(\pi(x),\Phi,\Psi)\ge0$ if $\pi(x)$ is consistent.
\item $D(\pi(x),\Phi,\Psi)\ge n+1$ if there is a sequence $(a_i:i<\omega)$ such
that $\models\Psi(\bar a)$ for any $k$-tuple $\bar a\subset(a_i:i<\omega)$, and 
$D(\pi(x)\land\Phi(x,a_i),\Phi,\Psi)\ge n$ for all $i<\omega$.
\end{itemize}
If $D(\pi,\Phi,\Psi)\ge n$ for all $n<\omega$, we put $D(\pi,\Phi,\Psi)=\infty$.
\end{enumerate}
An {\em inconsistency witness} is a $k$-inconsistency witness, for some $k<\omega$.
\edefn
\bem\label{indisc} Note that $D(\pi(x,a),\Phi,\Psi)\ge n$ is a closed condition 
on~$a$, and
$D(\tp(x/a),\Phi,\Psi)\ge n$ is a closed condition on $x$ over $a$. By
compactness and Ramsey's theorem, we may require $(a_i:i<\omega)$ to be
indiscernible in Definition \ref{Ddef}\,(\ref{Ddefdef}).\ebem
\lmm\label{Dinfty} Let $\Psi$ be an inconsistency witness for $\Phi$, and 
$\pi$ a partial type over $A$. Then
$D(\pi,\Phi,\Psi)$ is infinite if and only if for every linear order $I$
there are elements $(b_i,a_i^j:i\in I, j<\omega)$ such that $\models\Psi(\bar
a)$ for all $\bar a\subset(a_i^j:j<\omega)$ of the right length,
$b_i\models\pi\land\bigwedge_{k\le i}\Phi(x,a_k^0)$, and $(a_i^j:j<\omega)$ is
indiscernible over $A\cup\{b_ka_k^0:k<i\}$, for all $i\in I$. Moreover, we 
may require $(b_ia_i^0:i\in I)$ to be indiscernible.\elmm
\bew If the condition is satisfied, we can take $I=\omega$. Then for all
$n\in\omega$ the partial type $\pi\land\bigwedge_{i\le n}\Phi(x,a_i^0)$ is
satisfied by $b_n$ and hence non-empty. So
$$\begin{aligned}D(\pi,\Phi,\Psi)&>D(\pi\land\Phi(x,a_0^0),\Phi,
\Psi)>D(\pi\land\Phi(x,a_0^0)\land\Phi(x,a_1^0),\Phi,\Psi)\\
&>\cdots>D(\pi\land\Phi(x,a_0^0)\land\cdots\land\Phi(x,a_n^0),\Phi,\Psi)\ge0.
\end{aligned}$$
Hence $D(\pi,\Phi,\Phi)>n$ for all $n<\omega$, and
$D(\pi,\Phi,\Psi)=\infty$.

For the converse, by compactness it is sufficient to consider finite $I$. We
show by induction that if $D(\pi,\Phi,\Psi)\ge n$, then the condition is
satisfied for $I$ of size $n$.
For $n=0$ there is nothing to show. Suppose $D(\pi,\Phi,\Psi)\ge n+1$. Then
by definition there is a sequence $(a_0^j:j<\omega)$ whose 
subsequences satisfy $\Psi$, and such that 
$D(\pi\land\Phi(x,a_0),\Phi,\Psi)\ge n$ for all $j<\omega$. By Remark 
\ref{indisc} we may assume that $(a_0^j:j<\omega)$ is indiscernible over $A$. 
Choose $b_0\models\pi\land\Phi(x,a_0^0)$.
By inductive hypothesis for the partial type $\pi\land\Phi(x,a_0^0)$ over 
$A\cup\{b_0,a_0^0\}$, there are $(b_i,a_i^j:1\le i\le n,j<\omega)$ such that 
$\models\Psi(\bar a)$ for all $\bar a\subset(a_i^j:j<\omega)$ of the right 
length, 
$$b_i\models\pi\land\Phi(x,a_0^0)\land\bigwedge_{1\le k\le i}\Phi(x,a_k^0),$$
and $(a_i^j:j<\omega)$ is indiscernible over 
$A\cup\{b_0,a_0^0\}\cup\{b_k,a_k^0:1\le k<i\}$, for all $1\le i\le n$, as 
required.

The final assertion follows by compactness and Ramsey's theorem.\qed
\defn Let $I$ be an ordered set. A sequence $I=(a_i:i\in I)$ is {\em independent
over $A$}, or {\em $A$-independent}, if $\tp(a_i/A,a_j:j<i)$ does not divide 
over $A$ for all $i\in I$. If $A\subseteq B$ and $p\in S(B)$, the sequence
$(a_i:i\in I)$ is a {\em Morley sequence in $p$ over $A$} if it is
$B$-indiscernible, $a_i\models p$ and
$\tp(a_i/B,a_j:j<i)$ does not divide over $A$ for all $i\in I$. If $A=B$, we
simply call it a {\em Morley sequence} in $p$.\edefn
\tats\label{presym}\cite[Corollary 3.2.5]{wa00} or \cite[Proposition 16.12]{C}
If $\tp(b/cd)$ does not divide over $d$ and $\tp(a/cbd)$ does not divide over
$bd$, then $\tp(ab/cd)$ does not divide over $d$.\etats 
For the rest of the section we fix a hyperdefinable set $X$ over $\emptyset$. 
We call a type $p(x)$ an {\em $X$-type} if it implies $x\in X$. 

The following 
theorem generalizes \cite[Theorem 2.4.7]{wa00} to the local hyperdefinable 
context. Note that in the classical development the forking properties 
for hyperimaginaries are deduced from the corresponding properties for 
representatives. Here we cannot do this, as the ambient theory may well not be 
simple. So we have to work with hyperimaginaries in $X$ throughout.
\satz\label{simple} The following are equivalent:\begin{enumerate} 
\item {\em Symmetry} holds on $X$: For all $a,b,c\in X$, $\tp(a/bc)$ does not
divide over $b$ if and only if $\tp(c/ab)$ does not divide over $b$. 
\item {\em Transitivity} holds on $X$: If $a,b,c,d\in X$, then $\tp(a/bcd)$ does
not divide over $b$ if and only if $\tp(a/bc)$ does not divide over $b$ and
$\tp(a/bcd)$ does not divide over $bc$. 
\item {\em Local character} holds on $X$: There is $\kappa$ such that for all
countable $a\in X$ and $A\subset X$ there is $A_0\subseteq A$ with
$|A_0|\le\kappa$ such that $\tp(a/A)$ does not divide over $A_0$. In fact, we
can take $\kappa=2^{|T|}$. 
\item $D(.,\Phi,\Psi)<\infty$ for any partial $X$-type $\Phi(x,y)$ and
inconsistency witness $\Psi$ for $\Phi$. 
\item For any $A\subseteq B\subset X$, a partial $X$-type $\pi(x,B)$
does not divide over $A$ if and only if there is a Morley sequence $I$ in
$\tp(B/A)$ such that $\{\pi(x,B'):B'\in I\}$ is consistent. 
\end{enumerate}
If any of these conditions is satisfied, then for all $A\subseteq B\subset
X$ and $a\in X$ the type $\tp(a/B)$ does not divide over $A$ if and only if 
$$D(\tp(a/B),\Phi,\Psi)=D(\tp(a/A),\Phi,\Psi)$$ for all $(\Phi,\Psi)$. Moreover,
{\em Extension} holds on $X$:
For any partial $X$-type $\pi(x)$ over $B$, if $\pi$ does not divide over $A$
then it has a completion which does not divide over $A$.\esatz 
\bew $(1)\Rightarrow(2)$ Clearly, if
$\tp(a/bcd)$ does not divide over $b$, it
does not divide over $bc$ and $\tp(a/bc)$ does not divide over $b$. Conversely,
suppose that $\tp(a/bcd)$ does not divide over $bc$ and $\tp(a/bc)$ does not
divide over $b$. By symmetry, $\tp(d/abc)$ does not divide over $bc$ and
$\tp(c/ab)$ does not divide over $b$. By Fact \ref{presym} $\tp(cd/ab)$ does not
divide over $b$, so again by symmetry $\tp(a/bcd)$ does not divide over $b$. 

$(2)\Rightarrow(4)$ Suppose there is a partial $X$-type $\Phi$ and an
inconsistency witness $\Psi$ for $\Phi$ such that $D(x=x,\Phi,\Psi)=\infty$.
Put $I=\{\pm1,\pm(1+\frac1n):n>0\}$ and choose a sequence $(b_i,a_i^j:i\in 
I,j<\omega)$ as given by Lemma \ref{Dinfty}. Let $A^-=\{b_ia_i^0:i<-1\}$ and 
$A^+=\{b_ia_i^0:i>1\}$. Then $\tp(b_1/A^-A^+)$ does
not divide over $A^-$ and $\tp(b_1/A^-A^+a_{-1})$ does not divide over 
$A^-A^+$, since the
former is finitely satisfiable in $A^-$ and the latter in $A^+$. However,
$(a_{-1}^j:j<\omega)$ witnesses that $\Phi(x,a_{-1})$, and hence
$\tp(b_1/A^-A^+a_{-1})$, divides over $A^-$, contradicting transitivity.

$(4)\Rightarrow(3)$ Assume $(4)$. First, we note that for $A\subseteq
B\subset X$, if 
$$D(\tp(a/B),\Phi,\Psi)=D(\tp(a/A),\Phi,\Psi)$$
for all $(\Phi,\Psi)$, then $\tp(a/B)$ does not divide over $A$. This is
obvious, as if some $A$-indiscernible sequence $(B_i:i<\omega)$ in $\tp(B/A)$
witnesses dividing, we can take $\Phi(x,y)=\tp(a,B)$ and
$\Psi=\tp(B_1,B_2,\ldots,B_n)$ for $n<\omega$ sufficiently large. Then
$\Psi$ is an $n$-inconsistency witness (clearly, we may restrict to countable
$B$), and 
$$D(\tp(a/B),\Phi,\Psi)<D(\tp(a/A),\Phi,\Psi).$$
Given $\tp(a/A)$ it is hence enough to take $A_0\subseteq A$ big enough such
that
$$D(\tp(a/A),\Phi,\Psi)=D(\tp(a/A_0),\Phi,\Psi)$$
for all $(\Phi,\Psi)$. There are only $2^{|T|}$ such pairs, so we need at most
that many parameters. 

$(3)\Rightarrow(4)$ Suppose $D(x=x,\Phi,\Psi)=\infty$. Then for any cardinal
$\kappa$ we can find an indiscernible sequence
$(b_i,a_i^j:i\le\kappa^+,j<\omega)$ as in
Lemma \ref{Dinfty}. Since $\Phi(x,a_i^0)$ divides over $\{a_j^0:j<i\}$ for all 
$i\le\kappa^+$, the type $\tp(b_{\kappa^+}/a_i^0:i<\kappa^+)$ divides over any 
subset of its domain of cardinality $\le\kappa$. 

$(4)\Rightarrow(5)$. Assume (4). Given $a_E\in X$ and $A\subseteq B=b_E\subset
X$, for any pair $(\Phi,\Psi)$ and any formula $\p(y,b)$ we can adjoin either
$$\exists yz\,[xEy\land zEb\land\p(y,z)]\qquad\text{or}\qquad\exists 
yz\,[xEy\land zEb\land\neg\p(y,z)]$$ 
and preserve
$D(.,\Phi,\Psi)$-rank. By compactness we can thus complete $\tp(a_E/A)$ to an 
$E$-type $p$ over $B$ of the same $D(.,\Phi,\Psi)$-rank. In particular, no
$\Phi$-instance in $p$ divides over $A$ with $\Psi$ as inconsistency
witness. Coding finitely many pairs $(\Phi_i,\Psi_i:i<n)$ in a
single one, one obtains an extension $p$ such that no
$\Phi_i$-instance $\Psi_i$-divides for any $i<n$; by compactness we can do
this for all pairs $(\Phi,\Psi)$ simultaneously and obtain an
extension which does not divide over $A$. Take $B=X^\M\supset A$ for some 
sufficiently saturated model $\M$. Then a sequence
$(a_i:i<\omega)\subset B$ such that $a_i\models p\restriction_{(A,a_j:j<i)}$ is
a Morley sequence in $\tp(a/A)$.

This shows in particular that if $\pi(x,B)$ does not divide over $A$, then there
is a Morley sequence $I$ in $\tp(B/A)$ such that $\{\pi(x,B'):B'\in I\}$ is
consistent. 

Conversely, suppose that $\pi(x,B)$ divides over $A$, as witnessed by an
$A$-indiscernible sequence $(B_i:i<\omega)$ in $\tp(B/A)$
with $\bigcup_{i<\omega}\pi(x,B_i)$ inconsistent. Take any Morley
sequence $I$ in $\tp(B/A)$. By \cite[Corollary
2.2.8]{wa00} (which is shown there for real tuples, but transfers easily to
hyperimaginaries) we may assume that $B_i\widehat\ I$ is $A$-conjugate to $I$ 
for
all $i\in I$ and that $(B_i:i<\omega)$ is indiscernible over $AI$. If
$\bar\pi(x)=\bigcup_{B'\in I}\pi(x,B')$ were consistent, then $(B_i:i<\omega)$
would witness that 
$$D(\bar\pi(x)\land\pi(x,B_0),\pi(x,y),\Psi)<D(\bar\pi(x),\pi(x,y),\Psi)$$
for some inconsistency witness $\Psi$. But by $A$-conjugacy the two ranks must
be equal, a contradiction. 

$(5)\Rightarrow(1)$ Let us first show {\em Extension}. If $A\subseteq
B\subset X$ and $\pi(x,B)$ is a partial $X$-type which does not divide over
$A$, let $(B_i:i<\alpha)$ be a very long Morley sequence in $\tp(B/A)$.
Consider
any realization $a\models\bigwedge_{i<\alpha}\pi(x,B_i)$.
Since $\alpha$ is large, there is an infinite subset $J\subset\alpha$ such
that $\tp(B_i/aA)$ is constant for $i\in J$. Put $p(x)=\tp(a/AB)$, a completion of $\pi$. Then $(B_i:i\in J)$ witnesses that $p$ does not divide over $A$. 

Now given $a,b,c\in X$ such that $\tp(a/bc)$ does not divide over $b$, let
$B=X^\M\ni bc$ for some sufficiently saturated model $\M$, and $p$ an extension 
of $\tp(a/bc)$ to $B$ which does not divide over $b$. Choose a sequence
$(a_i:i<\omega)\subset B$ such that $a_i\models
p\restriction_{(bc,a_j:j<i)}$. This is a Morley sequence in $\tp(a/bc)$ over
$b$. Then $(a_i:i<\omega)$ is a Morley sequence in $\tp(a/b)$, and
$a_i\models\tp(a/bc)$ for all $i<\omega$. Hence $\tp(c/ba)$ does not divide over
$b$, and symmetry holds. 

Finally we show that if $(1)-(5)$ hold and $\tp(a/B)$ does not divide over $A$
for $A\subseteq B\subset X$ and $a\in X$,
then $D(\tp(a/A),\Phi,\Psi)\ge n$ implies $D(\tp(a/B),\Phi,\Psi)\ge n$ for all
$(\Phi,\Psi)$. For $n=0$ this is obvious. So suppose $D(\tp(a/A),\Phi,\Psi)\ge
n+1$. Then there is $(d_i:i<\omega)$ indiscernible over $A$ such that $\bar
d\models\Psi$ for all $\bar d\subset(d_i:i<\omega)$ of the right length, and
$D(\tp(a/A)\land\Phi(x,d_i),\Phi,\Psi)\ge n$ for all $i<\omega$. Let $q$ be a
completion of $\tp(a/A)\land\Phi(x,d_0)$ with $D(q,\Phi,\Psi)\ge n$. Clearly,
we may assume $a\models q$, and that $\tp(d_0/aB)$ does not divide over $aA$. As
$\tp(a/B)$ does not divide over $A$, by symmetry and transitivity $\tp(ad_0/B)$
does not divide over $A$, and $\tp(a/d_0B)$ does not divide over $d_0A$. By
induction hypothesis, $$D(\tp(a/d_0A),\Phi,\Psi)\ge
n\quad\text{implies}\quad D(\tp(a/d_0B),\Phi,\Psi)\ge n.$$
As $\tp(d_0/B)$ does not divide over $A$ and $(d_i:i<\omega)$ is
$A$-indiscernible, we may assume that it remains indiscernible over $B$. But
then it witnesses
$$D(\tp(a/B),\Phi,\Psi)\ge D(\tp(a/d_0B),\Phi,\Psi)+1\ge n+1.\qed$$

\defn An $A$-hyperdefinable set $X$ is {\em simple (over $A$)} if it satisfies
any of the conditions of Theorem \ref{simple} when we adjoin $A$ to the
language. If $X$ is simple over $A$ and $a,b,c\in X$, we shall say that $a$ and
$c$ are independent over $Ab$, written $a\ind_{Ab}c$, if $\tp(a/Abc)$ does not
divide over $Ab$.\edefn
Note that we only allow tuples and parameters from $A\cup X$. If $X$ is stably
embedded, we can of course allow parameters from anywhere. It is immediate from the definition that if $X$ is simple over $A$ and $B\subset X$, then $X$ is simple over $AB$.
\bem If $X$ is merely hyperdefinable, it may be simple although no definable or
even type-definable imaginary set in the ambient structure is simple.\ebem
If $X$ is simple, it is now standard to extend the notions of dividing and
independence to hyperimaginaries in $X_A^{heq}$. Moreover, we can
develop basic simplicity theory (canonical bases, the independence theorem,
stratified ranks, generic types, stabilizers, see \cite{C,wa00})
within $X_A^{heq}$, replacing models $\M$ by subsets $X_A^{heq}\cap\M^{heq}$.
\satzli\label{simint} Let $X$ and $Y$ be orthogonal $\emptyset$-hyperdefinable 
sets such that
$X$ is simple over $\emptyset$. If $A\subset Y$ is a set of parameters, then
$X$ is simple over $A$, and over $\bdd_Y(A)$. In particular, let $Z$ be a set
hyperdefinable over some parameters $A\subset X\cup Y$. If $Z$ is
$X$-internal within $X\cup Y$, then $Z$ is simple over $A$; if
$Z\subseteq(X\cup Y)^{heq}$ is almost $X$-internal, then $Z$ is simple
as well.\esatzli
\bew Simplicity over $A$ is obvious from orthogonality;
simplicity over $\bdd_Y(A)$ follows. Now if $Z$ is $X$-internal within $X\cup
Y$, then $Z\subset X_{\bdd_Y(A)}^{heq}$ by Proposition \ref{anaint}, and must be simple as well; if $Z\subseteq(X\cup 
Y)^{heq}$ is almost $X$-internal within $X\cup Y$, it is $X$-internal within
$\bdd_Y(A)$ by Corollary \ref{aint=>int}.\qed

\section{A hyperdefinable version of Schlichting's Theorem}\label{secfive}
Recall that Schlichting's Theorem \cite{Sch80}, generalized by Bergman and Lenstra \cite{BL89}, states that if $\H$ is a family of uniformly commensurable subgroups of a group $G$, then there is a subgroup $N$ commensurable with all groups in $\H$ (in fact a finite extension of a finite intersections of groups in $\H$) which is invariant under all automorphisms of $G$ which fix $\H$ setwise. Here two subgroups $H$ and $K$ are commensurable if their intersection has finite index in both $H$ and in $K$; uniformly commensurable means that there is a finite bound on these indices as $H$ and $K$ vary inside $\H$.

We shall call two hyperdefinable subgroups $G$ and $H$ {\em commensurable} if the index of their intersection in both $G$ and in $H$ is bounded, i.e.\ less than the cardinality $\kappa$ of the monster model. If $G$ is a hyperdefinable group, a {\em hyperdefinable family of subgroups} is a family $\H=\{H_a:a\models\pi\}$ for some partial types $\pi(y)$ and $\Phi(x,y)$ such that $H_a=\{x\in G:\models\Phi(x,a)\}$ is a subgroup of $G$ for any $a\models\pi$.
Note that a hyperdefinable family of commensurable subgroups is automatically uniformly commensurable by compactness, i.e.\ the index of the intersection $H\cap H^*$ in $H$ (and by symmetry in $H^*$) is bounded independently of the choice of $H,H^*\in\H$.

For hyperdefinable families of commensurable groups in a simple theory a version of Schlichting's Theorem has been shown in \cite[Theorem 4.5.13]{wa00}, generalizing a result of Hrushovski for theories of finite and definable S1-rank.
Here we shall show it for hyperdefinable families of commensurable subgroups in any theory.

\satz\label{schlichting} Let $G$ be a hyperdefinable group, $\H$ a hyperdefinable family of commensurable subgroups, and $\Gamma$ a hyperdefinable group of automorphisms of $G$ stabilizing $\H$ setwise. Then there is a $\Gamma$-invariant hyperdefinable subgroup $N$ commensurable with any group in $\H$; moreover $N$ is invariant under any model-theoretic automorphism stabilising $\H$.\esatz
\bew Suppose $\H=\{\Phi(x,a):a\models\pi\}$; clearly we may assume that $\Phi$ is closed under finite conjunctions. We put $H_a=\{g\in G:\models\Phi(g,a)\}$. As $\H$ is $\Gamma$-invariant, we have $\gamma H\in\H$ for any $H\in\H$.

Enumerate $\Phi=\{\phi_i:i<\alpha\}$ for some ordinal $\alpha$; for $i<\alpha$ put
$$\psi_i(x,x',y,\zeta)=\exists z\,[z\in\zeta^{-1}(x^{-1}x')\land\neg\phi_i(z,y)].$$
(Here $\zeta$ is a variable for elements from $\Gamma$, acting on elements from $G$. Thus $\zeta^{-1}(x^{-1}x')$ is an element of $G$ depending on $x$, $x'$ and $\zeta$; as it is hyperimaginary, it corresponds to a class of real tuples, and we demand that $z$ be one of them.)
Clearly $\psi_i$ is (equivalent to) a partial type. Consider a hyperdefinable subgroup $K$ of $G$. A {\em complete $\psi_i(x,x',a,\gamma)$-graph of size $n$ in $K$} is a set of elements $\{h_j:j<n\}$ of $K$ such that for any $j\not=j'$ one has $\models \psi_i(h_j,h_j',a,\gamma)$. The existence of such a graph implies in particular that the index of $K\cap\gamma H_a$ in $K$ is at least $n$.

As $\H$ is a hyperdefinable family of commensurable subgroups, by compactness for every $i<\alpha$ there is an integer $n_i$ such that for any $H\in\H$ and $(\gamma,a)\models\Gamma\times\pi$ there is no complete $\psi_i(x,x',a,\gamma)$-graph in $H$ of size $n_i$. If $K$ is a bounded intersection of groups in $\H$, we put $i(K,a,\gamma)=(k_i:i<\alpha)$, where $k_i\le n_i$ is the size of a maximal complete $\psi_i(x,x',a,\gamma)$-graph in $K$, and call this the {\em index} of $(a,\gamma)$ in $K$.  We order the set of indices lexicographically. 

Clearly, for $K_0\le K_1$ we have $i(K_0,a,\gamma)\le i(K_1,a,\gamma)$; by compactness equality holds if and only if $K_0\gamma H_a=K_1\gamma H_a$: if 
$i(K_0,a,\gamma)=i(K_1,a,\gamma)$, then for any $i<\alpha$ let $(g_j^i:j<k_i)$ be a maximal complete $\psi_i(x,x',a,\gamma)$-graph in $K_0$. By equality of the index, this is also a maximal complete graph in $K_1$. Then for any $g\in K_1$  and $i<\alpha$ there is some $g_j^i\in K_0$ with $\models\phi_i(\gamma^{-1}(g^{-1}g_j^i),a)$. By compactness, as $\Phi$ is closed under finite conjunctions, there is $g'\in K_0$ such that $\models\phi_i(\gamma^{-1}(g^{-1}g'),a)$ for all $i<\alpha$, that is $g^{-1}g'\in\gamma H_a$ and $g\in K_0\gamma H_a$.

By compactness, for every bounded intersection $K$ of groups in $\H$ there is some maximal index $i$ such that for some $(\gamma,a)\models\Gamma\times\pi$ we have $i=i(K,a,\gamma)$; call this the index $i(K)$ of $K$. Since for $\gamma'\in\Gamma$ and $(\gamma,a)\models\Gamma\times\pi$ we have
$$i(\gamma' K,a,\gamma'\gamma)=i(K,a,\gamma),$$
we obtain $i(K)=i(\gamma K)$. As the set of indices is bounded and $i(K_0)\le i(K_1)$ for $K_0\le K_1$, there is some bounded intersection $K$ of groups in $\H$ such that $i(K)$ is minimal possible, say $i_0$. We shall call $K$ {\em strong} if $i(K)=i_0$.

If $K$ is strong, we put 
$$\H(K)=\{H\in\H:\exists\,(\gamma,a)\models\Gamma\times\pi\ [H=\gamma H_a\land i(K,a,\gamma)=i_0]\}.$$
Then $\gamma' K$ is also strong for any $\gamma'\in\Gamma$, and
$$\H(\gamma' K)=\{\gamma'\gamma H_a:\gamma H_a\in\H(K)\}=\gamma'(\H(K)).$$
Now for $H\in\H(K)$ the set $\bigcap_{g\in K}(KH)^g$ is a subgroup of $G$ containing $K$; it is hyperdefinable, as we only have to conjugate by a set of representatives of $K/H$, which is bounded. Then 
$$N(K)=\bigcap_{H\in\H(K)}\bigcap_{g\in K}(KH)^g$$
is a subgroup of $G$ containing $K$; it is hyperdefinable as it contains $K$ and must have bounded index in $\bigcap_{g\in K}(KH)^g$ for any $H\in\H(K)$, so is a bounded intersection.

If $K_1$ is strong and $K_0\le K_1$, then $K_0$ is again strong and  $\H(K_0)\subseteq\H(K_1)$. Moreover $K_0H=K_1H$ for any $H\in\H(K_0)$, whence
$$\bigcap_{g\in K_0}(K_0H)^g=\bigcap_{g\in K_1}(K_1H)^g,$$
and 
$$K_1\le N(K_1)\le N(K_0)\le K_0H=K_1H.$$
It follows that there is an absolute bound on the index $|N(K):K_1|$, independent of the choice of strong $K$. As for a bounded family $(K_i:i\in I)$ of strong subgroups the intersection $\bigcap_{i\in I}K_i$ is again strong, there is some strong $K$ such that $N=N(K)$ is maximal possible. Then $N$ is hyperdefinable, commensurable with all groups in $\H$, and invariant under $\Gamma$ and all model-theoretic automorphisms stabilizing $\H$ setwise.\qed
\kor\label{normal} Let $G$ be a hyperdefinable group, and $H$ a subgroup commensurable with all its $G$-conjugates. Then there is a normal hyperdefinable subgroup $N$ commensurable with $H$.\ekor
\bew We apply Theorem \ref{schlichting} to the family $\H$ of $G$-conjugates of $H$, with the action of $\Gamma=G$ by conjugation.\qed

\section{Groups interpretable in orthogonal sets}
Recall that two hyperdefinable subgroups $H_1$ and $H_2$ of some group $G$ are
{\em commensurable} if $H_1\cap H_2$ has bounded index both in $H_1$ and in
$H_2$.
\satz\label{quotient} Suppose $X$ and $Y$ are orthogonal
$\emptyset$-hyperdefinable sets in a structure $\M$, and $G$ is an
$\emptyset$-hyperdefinable group in $(X\cup Y)^{heq}$. If $X$ is simple over
$\emptyset$, there is an $\emptyset$-hyperdefinable normal $X$-internal subgroup
$N$ of $G$ such that the quotient $G/N$ is $Y$-internal. $N$ is unique up to
commensurability.\esatz 
\bew Let us first show uniqueness: If $N'$ is a second such group, then
$N/(N\cap N')$ and $N'/(N\cap N')$ are $X$-internal and $Y$-internal, and hence
bounded by orthogonality of $X$ and $Y$. Thus $N$ and $N'$ are commensurable. 

By Theorem \ref{wehi} every element $g\in G$ is of the form $(g_X,g_Y)_E$ for
some $g_X\in X^{heq}$ and $g_Y\in Y^{heq}$, both bounded over $g$, and some
type-definable equivalence relation $E$  with bounded classes, depending on
$\tp(g)$. Hence $\tp(g/g_Y)$ is $X$-internal and $\tp(g/g_X)$ is
$Y$-internal. Now if $h=(h_X,h_Y)_E$ and $gh=((gh)_X,(gh)_Y)_E$, then
$$(gh)_X\in\bdd_X(g_X,g_Y,h_X,h_Y),$$
whence $(gh)_X\in\bdd_X(g_X,h_X)$ by Corollary
\ref{dclbdd}. Similarly $(gh)_Y\in\bdd_Y(g_Y,h_Y)$. 

Now $X$ is simple, as is $\tp(g,h/g_Y,h_Y)$ for any $g,h\in G$ by 
Proposition \ref{simint}. Hence we can
consider $g,h\in G$ such that $g\ind_{g_Y,h_Y}h$. Then for any stratified local
rank~$D$ 
\begin{equation}\tag{$\dag$}\begin{aligned}D(gh/(gh)_Y)&\ge 
D(gh/(gh)_Y,g_Y,h_Y,g)=D(h/g_Y,h_Y,g)\\
&=D(h/g_Y,h_Y)=D(h/h_Y),\end{aligned}\end{equation}
where the first equality holds as $(gh)_Y\in\bdd_Y(g_Y,h_Y)$, and the last 
equality follows from $h\ind_{h_Y}g_Y$ by orthogonality (Example 
\ref{ind}). Similarly
$$D(gh/(gh)_Y)\ge D(g/g_Y).$$
Now suppose $G$ is a subset of $(X^m\times Y^n)/E$, where $m,n$ are at
most countable. Then if $g=(\bar x_g,\bar y_g)_E\in G$, we have $g\ind_{g_Y}\bar
y_g$ by $X$-internality of $\tp(g/g_Y)$ and orthogonality, and $g_Y\in\bdd(\bar 
y_g)$,
whence $D(g/g_Y)=D(g/\bar y_g)$. By compactness, there is a $G$-type
$p((\bar x,\bar y)_E)$ implying that $D((\bar x,\bar y)_E/\bar y)$ is
maximal for all local stratified ranks. But if $g,h\models p$ with
$g\ind_{g_Y,h_Y}h$, then we must have equality in $(\dag)$. Therefore $g$, $h$ and $gh$ are pairwise independent over $g_Y,h_Y,(gh)_Y$. Put $A=(g_Y,h_Y,(gh)_Y)$. Then $X'=\tp(g/A)$ is $X$-internal and simple, as is $X'X'^{-1}$. We may therefore define
$$S_0=\{g\in G:\exists x\, (x\equiv^{lstp}_Agx\equiv^{lstp}a\land x\ind_Ag\land gx\ind_Ag\}\subseteq X'X'^{-1}$$
and the stabilizer $S=\St(g/A)=S_0^2$, an $X$-internal hyperdefinable subgroup of $G$.

Now \cite[Lemme 1.2]{bmpw} (see \cite[Remarque 1.3]{bmpw} for the extension from the stable to the simple context) states that whenever $g$, $h$ and $gh$ are pairwise independent over $A$, then $g$ is generic in the coset $S g$, and this coset is hyperdefinable over $\bdd(A)$. 

By orthogonality, $g\ind_{g_Y}h_Y,(gh)_Y$. This implies in particular 
$$D(S)=D(Sg)=D(g/g_Y,h_Y,(gh)_Y)=D(g/g_Y)=D(p).$$

Suppose that $S$ is not commensurable with $S^h$ for some
$h\in G$. Then  $SS^h$ is still $X$-internal, with $$D(SS^h)\ge D(S)=D(p)$$
for every stratified local rank $D$, and for at least one such rank $D_0$ we
have $D_0(SS^h)>D_0(p)$. Choose $g'\in SS^h$ with $D_0(g'/h)=D_0(SS^h)$. By
pre-multiplying with a generic element of $S$ and post-multiplying with a
generic element of $S^h$, the inequality $(\dag)$ implies that we may
assume $D(g'/h)\ge D(p)$ for every stratified local rank $D$. However,
$\tp(g'/h)$ is $X$-internal, so $g'_Y\in\bdd_Y(g')$ implies $g'_Y\in\bdd_Y(h)$ 
by Lemma \ref{dclbdd} and Proposition \ref{anaint}. Thus 
$$D(g'/g'_Y)\ge D(g'/h)\ge D(p)\quad\text{and}\quad D_0(g'/g'_Y)\ge
D_0(g'/h)>D_0(p),$$
contradicting our choice of $p$. Hence $S$ is commensurable with all its
conjugates. By Corollary \ref{normal} there is a hyperdefinable normal subgroup $N$ of $G$ commensurable with $S$. So $N$ is $X$-internal, and $D(N)=D(S)=D(p)$ for all stratified local ranks $D$.
Now the same proof, with $NZ$ instead of $SS^h$, shows that if $Z$ is an
$X$-internal hyperdefinable subset of $G$, then $Z$ is covered by boundedly many
cosets of $N$. In particular, for any $g'\in G$ the type $\tp(g'/g'_Y)$ is
covered by boundedly many cosets of $N$. But then $g'N\in\bdd(g'_Y)$, and $G/N$
is almost $Y$-internal, whence $Y$-internal by Corollary \ref{aint=>int}.\qed 

\kor\label{BM} Suppose $X$ and $Y$ are orthogonal type-definable sets over
$\emptyset$ in 
a structure $\M$, and $G$ is a type-interpretable group over $\emptyset$ in 
$(X\cup Y)^{eq}$. If $X$ is simple over
$\emptyset$, there is a normal $X$-internal subgroup $N$ of $G$ 
type-interpretable over $\emptyset$, such that the quotient $G/N$ is 
$Y$-internal. $N$ is unique up to commensurability. If $X$ is definable and 
supersimple, then we can take $N$ relatively interpretable.\ekor
\bew The first part is obvious from Theorem \ref{quotient}, as a 
hyperdefinable subgroup of a type-interpretable group is again 
type-interpretable.

If $X$ is definable and supersimple, $N$ must be contained in a
definable $X$-internal set $\bar X$ by \cite[Lemma 3.4.17]{wa00}; note
that $\bar X$ will also be supersimple. So $N$ is the intersection of definable
supergroups by \cite[Theorem 5.5.4]{wa00}, one of which, say $N_0$, must be
contained in $\bar X$ by compactness. Then $N_0$ is 
$X$-internal. As above, $N_0$ must be commensurable with all its 
$G$-conjugates; moreover, commensurability is uniform by compactness (or
\cite[Lemma 4.2.6]{wa00}). By \cite[Theorem 4.2.4]{wa00} there is a
relatively interpretable normal subgroup $\bar N$ commensurable with $N_0$. So
$\bar N$ is $X$-internal, and $G/\bar N$ is $Y$-internal.\qed
\begin{frage} If $X$ and $Y$ are orthogonal type-definable sets (or even 
definable sets), $X$ is simple and $G$ is a
relatively definable group in $(X\cup Y)^{eq}$, can we find a relatively
definable normal $X$-internal subgroup $N$ such that $G/N$ is
$Y$-internal?\end{frage} 

\begin{frage} What can we say if neither $X$ nor $Y$ is simple? Is it true that
in every hyperdefinable subgroup of $((X\cup Y)^{heq}$ there is a maximal normal
hyperdefinable $X$-internal subgroup $N_X$, a maximal normal hyperdefinable
$Y$-internal subgroup $N_Y$, an $X$-internal hyperdefinable local group $G_X$, a
$Y$-internal hyperdefinable local group $G_Y$ and a hyperdefinable
locally bounded equivalence relation $E$ on $G_X\times G_Y$ such that
$G/(N_XN_Y)$ is isogenous, or even isomorphic, to $(G_X\times
G_Y)/E$~?\end{frage}

\end{document}